\documentclass[12pt]{amsart}
\usepackage{amssymb}
\usepackage[all]{xy}

\newcommand{\issuenumber}{14}
\newcommand{\issuemonth}{September}
\newcommand{\issueyear}{2005}

\newtheorem{thm}{Theorem}[section]
\newtheorem{prob}[thm]{Problem}

\newtheorem{issue}{Issue}

\theoremstyle{definition}

\theoremstyle{remark}

\newcommand{\Fin}{{[\N]^{<\aleph_0}}}

\newcommand{\ed}{\end{thebibliography}\general\end{document}}

\newcommand{\x}{\times}

\newcommand{\Cantor}{{{}^\N\{0,1\}}}

\newcommand{\roth}{P_\infty(\N)}

\newcommand{\fd}{\mathfrak{d}}
\newcommand{\fp}{\mathfrak{p}}

\newcommand{\NON}{{\mathsf   {NON}}}
\newcommand{\COF}{{\mathsf   {COF}}}

\newcommand{\M}{\mathcal{M}}

\newcommand{\cov}{\mathsf{cov}}

\renewcommand{\c}{\mathfrak{c}}

\newcommand{\R}{\mathbb{R}}

\newcommand{\fo}{\mathfrak{od}}

\renewcommand{\b}{\mathfrak{b}}

\renewcommand{\split}{\mathsf{Split}}
\newcommand{\bq}{\begin{quote}}
\newcommand{\eq}{\end{quote}}
\renewcommand{\O}{\mathcal{O}}
\newcommand{\B}{\mathcal{B}}
\newcommand{\BG}{\B_\Gamma}

\newcommand{\sone}{\mathsf{S}_1}    \newcommand{\sfin}{\mathsf{S}_{fin}}
\newcommand{\ufin}{\mathsf{U}_{fin}}

\newcommand{\nin}{\not\in}
\newcommand{\cF}{\mathcal{F}}

\newcommand{\fu}{\mathfrak{u}}

\newcommand{\NN}{{{}^{\naturals}\naturals}}
\newcommand{\naturals}{{\mathbb N}}
\newcommand{\N}{\naturals}
\newcommand{\sm}{\setminus}
\newcommand{\sbst}{\subseteq}
\newcommand{\by}[2]{\par\hfill\emph{#1}, #2}
\newcommand{\nby}[1]{\par\hfill\emph{#1}}
\newcommand{\Tau}{\mathrm{T}}
\newcommand{\CE}{\textsc{CE}}

\newcommand{\be}{\begin{enumerate}}
\newcommand{\ee}{\end{enumerate}}
\newcommand{\bi}{\begin{itemize}}
\newcommand{\ei}{\end{itemize}}

\newcommand{\general}{\small\vfill\par\noindent\hrulefill\par
\noindent\textbf{Previous issues.} The first issues of this
bulletin, which contain general information (first issue), basic
definitions, research announcements, and open problems (all
issues) are available online, on \arx{math.GN/$x$}, where $x$ is
\texttt{0301011}, \texttt{0302062}, \texttt{0303057},
\texttt{0304087}, \texttt{0305367}, \texttt{0312140},
\texttt{0401155}, \texttt{0403369}, \texttt{0406411},
\texttt{0409072}, \texttt{0412305}, \texttt{0503631},
and \texttt{0508563},
respectively, for issues number $1$ to $13$.\\[0.1cm]
\textbf{Contributions.}
Please submit your contributions (announcements, discussions, and open problems)
by e-mailing us. It is preferred to write them
in \LaTeX{}.
The authors are urged to use as standard notation as possible, or otherwise give
the definitions or a reference to where the notation is explained.
Contributions to this bulletin would not require any transfer of copyright,
and material presented here can be published elsewhere.\\[0.1cm]
\textbf{Subscription.}
To receive this bulletin (free) to your
e-mailbox, e-mail us:\\
{boaz.tsaban@weizmann.ac.il}
}

\newcommand{\nArxPaper}[5]{\subsection{#2}{#4}\par\hfill{\arx{#1}}\par\hfill\emph{#3}}

\newcommand{\nAMSPaper}[5]{\subsection{#2}{#4}\par\hfill{\texttt{#1}}\par\hfill\emph{#3}}

\newcommand{\arx}[1]{\texttt{http://arxiv.org/abs/#1}}
\newcommand{\url}[1]{\bq\texttt{#1}\eq}
\newcommand{\online}[1]{The paper is available online at \url{#1}}

\title[$\mathcal{SPM}$ Bulletin \textbf{\issuenumber} (\issuemonth{} \issueyear)]{%
$\mathcal{SPM}$ Bulletin\\[0.5cm]
Issue number \issuenumber: \issuemonth{} \issueyear{} \CE{}}

\begin{document}
\maketitle

\tableofcontents

\section{Editor's note}

The second SPM meeting promises to be a most interesting event.
It suffices to skim the list of participants.
Visit
\url{http://www.matematica.unile.it/ricerca/CSGT/workshop.htm}
and
\url{http://diamond.boisestate.edu/~spm/Lecce2/index.htm}
for details.
Those who have not yet submitted their abstract to this meeting
are advised to do so at their earliest convenience, at
\url{http://atlas-conferences.com/cgi-bin/abstract/submit/caqh-01}
Once the abstracts are submitted, it will be possible to see
them at
\url{http://at.yorku.ca/cgi-bin/amca/caqh-01}

\medskip

Contributions to the next issue are, as always, welcome.

\medskip

\by{Boaz Tsaban}{boaz.tsaban@weizmann.ac.il}

\hfill \texttt{http://www.cs.biu.ac.il/\~{}tsaban}

\section{Research announcements}

\nArxPaper{math.LO/0505471}
{New reals: Can live with them, can live without them}
{Martin Goldstern and Jakob Kellner}
{We give a simple proof of the preservation theorem for proper countable
support iterations known as ``tools-preservation'', ``Case A'' or ``first
preservation theorem'' in the literature. We do not assume that the forcings
add reals.

  Currently the paper contains drafts of the proofs only. We will fill in more
details soon.
}

\nArxPaper{math.LO/0506019}
{Uniform almost everywhere domination}
{Peter Cholak, Joseph Miller, and Noam Greenberg}
{We explore the interaction between Lebesgue measure and dominating functions.
We show, via both a priority construction and a forcing construction, that
there is a function of incomplete degree that dominates almost all degrees.
This answers a question of Dobrinen and Simpson, who showed that such functions
are related to the proof-theoretic strength of the regularity of Lebesgue
measure for $G_\delta$ sets. Our constructions essentially settle the reverse
mathematical classification of this principle.}

\nArxPaper{math.GN/0505516}
{Heredity of $\tau$-pseudocompactness}
{Jerry E.\ Vaughan}
{S.\ Garcia-Ferreira and H.\ Ohta gave a construction that was intended to
produce a tau-pseudocompact space, which has a regular-closed zero set $A$ and a
regular-closed $C$-embedded set $B$ such that neither $A$ nor $B$ is $\tau$-pseudocompact.
We show that although their sets $A, B$ are not regular-closed, there are at
least two ways to make their construction work to give the desired example.
}

\nArxPaper{math.LO/0505645}
{Understanding preservation theorems: $\omega^\omega$-bounding}
{Chaz Schlindwein}
{This is an expository note giving Shelah's proof of the preservation of
``proper + $\omega^\omega$-bounding'' (Theorem 1.12 of Chapter VI of his book
\emph{Proper and Improper Forcing}).
}

\nAMSPaper{http://www.ams.org/journal-getitem?pii=S0002-9947-05-03956-5}
{Classification problems in continuum theory}
{Riccardo Camerlo, Udayan B.\ Darji, and Alberto Marcone}
{We study several natural classes and relations occurring in
continuum theory from the viewpoint of descriptive set theory and
infinite combinatorics. We provide useful characterizations for
the relation of likeness among dendrites and show that it is a bqo
with countably many equivalence classes. For dendrites with
finitely many branch points the homeomorphism and
quasi-homeomorphism classes coincide, and the minimal quasi-homeomorphism
classes among dendrites with infinitely many branch points are
identified. In contrast, we prove that the
homeomorphism relation between dendrites is $S_\infty$-universal. It is
shown that the classes of trees and graphs are both
$\mathrm{D}_{2}({{\boldsymbol
\Sigma_{3}^{0}}})$-complete, the class of dendrites is
${{\boldsymbol\Pi_{3}^{0}}}$-complete, and the
class of all continua homeomorphic to a graph or dendrite with
finitely many branch points is
${{\boldsymbol\Pi_{3}^{0}}}$-complete. We also show that
if $G$ is a nondegenerate finitely triangulable continuum, then
the class of $G$-like continua is
${\boldsymbol\Pi_{2}^{0}}$-complete.
}

\nArxPaper{math.GN/0506379}
{Residuality of families of $F_\sigma$ sets}
{Shingo Saito}
{We prove that two natural definitions of residuality of families of $F_\sigma$
sets are equivalent.
We make use of the Banach-Mazur game in the proof.}

\nAMSPaper{http://www.ams.org/journal-getitem?pii=S0002-9947-05-04034-1}
{First countable, countably compact spaces and the continuum hypothesis}
{Todd Eisworth and Peter Nyikos}
{We build a model of ZFC+CH in which every first countable, countably compact
space is either compact or contains a homeomorphic copy of $\omega_1$ with the
order topology. The majority of the paper consists of developing forcing
technology that allows us to conclude that our iteration adds no reals. Our
results generalize Saharon Shelah's
iteration theorems appearing in Chapters V and VIII of
{\it Proper and improper forcing} (1998), as well as Eisworth and Roitman's
(1999) iteration theorem.
We close the paper with a ZFC example (constructed using Shelah's
club--guessing sequences) that shows similar results do not hold for closed
pre--images of $\omega_2$.
}

\subsection{Game approach to universally Kuratowski-Ulam spaces}
We consider a variant of the open-open game. A topological
characterization of I-favorable spaces is used to show that
the hyperspace over I-favorable space is I-favorable. The main result
is that every I-favorable space is universally Kuratowski-Ulam.\footnote{A
space $Y$ is universally Kuratowski-Ulam if the Kuratowski-Ulam Theorem
holds in $X\x Y$ for every space $X$.}
\nby{Andrzej Kucharski and Szymon Plewik}

\nArxPaper{math.GN/0507062}
{Small Valdivia compact spaces}
{Wies{\l}aw Kubi\'s and Henryk Michalewski}
{We prove a preservation theorem for the class of Valdivia compact
spaces, which involves inverse sequences of ``simple''
retractions. Consequently, a compact space of weight
$\leqslant\aleph_1$ is Valdivia compact iff it is the limit of an
inverse sequence of metric compacta whose bonding maps are
retractions. As a corollary, we show that the class of Valdivia
compacta of weight $\leqslant\aleph_1$ is preserved both under
retractions and under open $0$-dimensional images. Finally, we
characterize the class of all Valdivia compacta in the language of
category theory, which implies that this class is preserved under
all continuous weight preserving functors.
}

\nAMSPaper{http://www.ams.org/journal-getitem?pii=S0002-9947-05-04000-6}
{Canonical forms of Borel functions on the Milliken space}
{Olaf Klein and Otmar Spinas}
{The goal of this paper is to canonize Borel measurable mappings
$\Delta\colon\Omega^\omega\to\mathbb{R}$, where $\Omega^\omega$ is the Milliken
space, i.e., the space of all increasing infinite sequences of pairwise
disjoint nonempty finite sets of $\omega$. This main result is a common
generalization of a theorem of Taylor and a theorem of Pr\"omel and Voigt.}

\nAMSPaper{http://www.ams.org/journal-getitem?pii=S0002-9947-05-04005-5}
{Complete analytic equivalence relations}
{Alain Louveau and Christian Rosendal}
{We prove that various concrete analytic equivalence relations
arising in model theory or analysis are complete, i.e.\ maximum in
the Borel reducibility ordering. The proofs use some general
results concerning the wider class of analytic quasi-orders.
}

\subsection{The Number of Near-Coherence Classes of Ultrafilters is Either Finite or $2^\c$}
We prove that the number of near-coherence classes of non-principal ultrafilters on the natural
numbers is either finite or $2^\c$.
Moreover, in the latter case the Stone-\v{C}ech compactification $\beta\N$ of $\N$ contains a
closed subset $C$ consisting of $2^\c$ pairwise non-nearly-coherent ultrafilters.
We obtain some additional information about such closed sets
under certain assumptions involving the cardinal characteristics $\fu$ and
$\fd$.

Applying our main result to the Stone-\v{C}ech remainder $\beta\R_+\sm\R_+$
of the half-line $\R_+ = [0,\infty)$ we obtain that the number of composants
of $\beta\R_+\sm\R_+$ is either finite or $2^\c$.
\nby{Taras Banakh and Andreas Blass}

\nAMSPaper{http://www.ams.org/journal-getitem?pii=S0002-9947-05-04074-2}
{The complexity of recursion theoretic games}
{Martin Kummer}
{
We show that some natural games introduced by Lachlan in 1970 as a model
of recursion theoretic constructions are undecidable, contrary to what was
previously conjectured. Several consequences are pointed out; for instance, the set of all
$\Pi_2$-sentences that are uniformly valid in the lattice of recursively enumerable sets is undecidable.
Furthermore we show that these games are equivalent to natural subclasses of effectively presented Borel games.
}

\nArxPaper{math.GN/0509097}
{Cosmic dimension}
{Alan Dow and Klaas Pieter hart}
{Martin's Axiom for $\sigma$-centered partial orders implies that there is a
cosmic space with non-coinciding dimensions.}

\nArxPaper{math.GN/0509099}
{There is no categorical metric continuum}
{Klaas Pieter Hart}
{We show there is no categorical metric continuum. This means that for every
metric continuum X there is another metric continuum Y such that X and Y have
(countable) elementarily equivalent bases but X and Y are not homeomorphic. As
an application we show that the chainability of the pseudoarc is not a
first-order property of its lattice of closed sets.
}

\section{On a conjecture and a problem of Hurewicz}\label{SFH}

A set of reals $X$ has \emph{Menger's property} (1924) if
no continuous image of $X$ in $\NN$ is cofinal with respect to $\le^*$.
It has the formally stronger \emph{Hurewicz' property} (1925) if
every continuous image of $X$ in $\NN$ is bounded.
$\sigma$-compact\-ness implies Hurewicz' property, which implies Menger's.
Both Menger and Hurewicz conjectured that their property characterizes $\sigma$-compactness,
and for a long time only consistent counter-examples were known.
Hurewicz (1927) also posed the problem whether
there is $X\sbst\R$ which is Hurewicz but not Menger.
The problem was raised again by Bukovsk\'y and Hale\v{s} (2003).

Fremlin-Miller (1988) and then Just-Miller-Scheepers-Szeptycki (1996)
gave a dichotomic existential argument refuting the
Conjectures in ZFC. Using the Michael topological technique, Chaber-Pol (2002)
improved the dichotomic argument and essentially solved the Hurewicz
Problem, alas in an existential manner.

Barotszy\'nski-Tsaban (2002) gave two explicit counter-examples to the conjectures
using two specialized constructions. Tsaban-Zdomsky (2005) generalize
both constructions and solve the Hure\-wicz Problem constructively by considering
scales with respect to semifilters (collections
of infinite subsets of $\N$ closed under almost supersets).
Working in $P(\N)$ (which is like $\Cantor$):
For each feeble semifilter $\cF$ and each $\cF$-scale $S$, all finite powers of $X=S\cup\Fin$
are Hurewicz and not $\sigma$-compact.
Viewed appropriately as a subset of $\R$, the field generated by $X$ is Hurewicz,
universally null, and universally meager.
The Hurewicz problem is solved by using
the semifilter $\cF=\roth$, and choosing the $\cF$-scale's points
such that (the enumerations of) their complements form an unbounded set.
To carry this, descriptive set theoretic properties of semifilters are used.
When $\fd$ is regular, subfields of $\R$ are constructed which are Menger but not Hurewicz.
This implies (in a dichotomic manner) that there always are such fields.

All results can be viewed as dealing with the Ramsey theory of open covers.

The work of Tsaban-Zdomsky, named \emph{Scales, fields, and a problem of Hurewicz},
is available at
\bq
\arx{math.GN/0507043}
\eq
\nby{Boaz Tsaban}

\section{Problem of the Issue}
The most interesting problem which arises from the work described in
Section \ref{SFH} is the following.

\begin{prob}
Does there exist (in ZFC) a set of reals $X$ of cardinality $\fd$ such that all
finite powers of $X$ have Menger's property $\ufin(\O,\O)$?
\end{prob}

We know that the answer is ``Yes'' when $\fd$ is regular.
If the answer is positive, then an explicit (non-dichotomic) construction
of such a set would be even more interesting.

\nby{Boaz Tsaban}

\section{Problems from earlier issues}

\begin{issue}
Is $\binom{\Omega}{\Gamma}=\binom{\Omega}{\Tau}$?
\end{issue}

\begin{issue}
Is $\ufin(\Gamma,\Omega)=\sfin(\Gamma,\Omega)$?
And if not, does $\ufin(\Gamma,\Gamma)$ imply
$\sfin(\Gamma,\Omega)$?
\end{issue}

\stepcounter{issue}

\begin{issue}
Does $\sone(\Omega,\Tau)$ imply $\ufin(\Gamma,\Gamma)$?
\end{issue}

\begin{issue}
Is $\fp=\fp^*$? (See the definition of $\fp^*$ in that issue.)
\end{issue}

\begin{issue}
Does there exist (in ZFC) an uncountable set satisfying $\sone(\BG,\B)$?
\end{issue}

\stepcounter{issue}

\begin{issue}
Does $X \nin \NON(\M)$ and $Y\nin\mathsf{D}$ imply that
$X\cup Y\nin \COF(\M)$?
\end{issue}

\begin{issue}
Is $\split(\Lambda,\Lambda)$ preserved under taking finite unions?
\end{issue}
\begin{proof}[Partial solution]
Consistently yes (Zdomsky). Is it ``No'' under CH?
\end{proof}

\begin{issue}
Is $\cov(\M)=\fo$? (See the definition of $\fo$ in that issue.)
\end{issue}

\begin{issue}
Does $\sone(\Gamma,\Gamma)$ always contain an element of cardinality $\b$?
\end{issue}

\begin{issue}
Could there be a Baire metric space $M$ of weight $\aleph_1$ and a partition
$\mathcal{U}$ of $M$ into $\aleph_1$ meager sets where for each ${\mathcal U}'\subset\mathcal U$,
$\bigcup {\mathcal U}'$ has the Baire property in $M$?
\end{issue}

\stepcounter{issue} 

\general\end{document}